\documentclass[12pt]{amsart}
\usepackage{amscd,amssymb,amsthm,verbatim}
\newcommand{\diam}{\operatorname{diam}}
\newcommand{\Lip}{\operatorname{Lip}}

\newcommand{\R}{{\mathbb R}}
\newcommand{\Ric}{\operatorname{Ric}}

\newcommand{\Z}{{\mathbb Z}}

\numberwithin{equation}{section}
\theoremstyle{plain}
\newtheorem{definition}[equation]{Definition}

\newtheorem{lemma}[equation]{Lemma}
\newtheorem{theorem}[equation]{Theorem}

\theoremstyle{definition}

\newtheorem{remark}[equation]{Remark}

\usepackage{graphicx}
\input{amssym.def}
\input{amssym}
\input{epsf}
\usepackage{a4wide}

\def\R{\mathbb R}

\def\Z{\mathbb Z}

\def\var{\varepsilon}

\def\dps{\displaystyle}

\def\2dr#1#2{\left. \frac{d^2}{d{#1}^2} \right |_{#2}}
\def\d2#1{\frac{d^2}{d{#1}^2}}

\def\LSI{{\rm LSI}}

\def\med{\medskip}
\def\sm{\smallskip}
\def\bul{$\bullet$ \ }

\def\begeq{\begin{equation}}
\def\endeq{\end{equation}}
\def\begar{\begin{eqnarray}}
\def\endar{\end{eqnarray}}
\def\begar*{\begin{eqnarray*}}
\def\endar*{\end{eqnarray*}}
\def\begal{\begin{align}}
\def\endal{\end{align}}
\def\begal*{\begin{align*}}
\def\endal*{\end{align*}}

\theoremstyle{definition}

\theoremstyle{remark}

\newtheorem*{Thm*}{Theorem}
\newtheorem*{Lem*}{Lemma}
\newtheorem*{Conj*}{Conjecture}
\newtheorem*{Cor*}{Corollary}
\newtheorem*{Def*}{Definition}
\newtheorem*{Prop*}{Proposition}
\newtheorem*{Exo*}{Exercise}
\newtheorem*{Exs*}{Examples}
\newtheorem*{Ex*}{Example}
\newtheorem*{Rk*}{Remark}
\newtheorem*{Rks*}{Remarks}

\begin{document}

\title[The Hamilton--Jacobi semigroup on length spaces and applications]
{Hamilton--Jacobi semigroup on length spaces and applications}

\author{John Lott}
\address{Department of Mathematics\\
University of Michigan\\
Ann Arbor, MI  48109-1109\\
USA} \email{lott@umich.edu}
\author{C\'edric Villani}
\address{UMPA\\ ENS Lyon\\
46 all\'ee d'Italie, 69364 Lyon Cedex 07\\
FRANCE} \email{cvillani@umpa.ens-lyon.fr}

\thanks{The research of the first author was 
supported by NSF grant DMS-0604829}
\date{March 28, 2007}

\begin{abstract}
We define a Hamilton--Jacobi semigroup acting on continuous functions 
on a compact length space. Following a strategy of Bobkov, Gentil
and Ledoux, we use some basic properties of the semigroup to study
geometric inequalities related to concentration of measure. Our main results
are that (1) a Talagrand inequality on a measured length space implies a global
Poincar\'e inequality and (2) if the space satisfies a doubling condition,
a local Poincar\'e inequality and a log Sobolev inequality  then it also
satisfies a Talagrand inequality.
\end{abstract}

\maketitle

Links between concentration of measure, log Sobolev inequalities,
Talagrand inequalities and Poincar\'e inequalities have been studied in the setting
of Riemannian manifolds~\cite{Toulouse (2000),Bobkov-Gentil-Ledoux 
(2001),Bobkov-Gotze (1999),Ledoux (1999), Ledoux (2001),Otto-Villani (2000)}.
The main result in the paper of Otto and Villani \cite{Otto-Villani (2000)} can be informally stated
as follows: on a Riemannian manifold, a log Sobolev inequality
implies a Talagrand inequality, which in turn implies a Poincar\'e
(or spectral gap) inequality, all of this being without any degradation
of the constants.

On the other hand, there has been intense recent activity
to develop a theory of Ricci curvature bounds, log Sobolev inequalities and
related inequalities in the more general setting of metric-measure
length spaces satisfying minimal regularity 
assumptions~\cite{Lott-Villani,Lott-Villani 2,von Renesse,Sturm1,Sturm2}. 

The goal of the present paper is to extend the main results
of~\cite{Otto-Villani (2000)} to this generalized framework, which can
be considered to be a natural degree of regularity for the problem.
To do so, we  adapt the strategy of  Bobkov-Gentil-Ledoux
\cite{Bobkov-Gentil-Ledoux (2001)},
based on the Hamilton--Jacobi semigroup. We also
establish the basic properties of the Hamilton--Jacobi semigroup
for general length spaces, which is of independent interest.

\section{Main results}

Basic information on length spaces
is in~\cite[Chapter~2]{Burago-Burago-Ivanov (2001)}.
For the sake of simplicity we work with  compact length spaces,
but the results remain valid for locally compact
complete separable length spaces.

Throughout this paper, $X$ will denote a compact length space,
equipped with a metric $d$ and a Borel reference probability
measure $\nu$. We use the following conventions:
\sm

- $\Lip(X)$ denotes the set of real-valued Lipschitz functions on $X$.
\sm

- Given $f \in C(X)$, we define the {\em gradient norm} of $f$ 
at a point $x \in X $ by
\begin{equation} \label{nablaf}
|\nabla f| (x) \: = \:
\limsup_{y\to x} \frac{|f(y)-f(x)|}{d(x,y)}.
\end{equation}
If $f \in \Lip(X)$ then $|\nabla f| \in L^\infty(X)$.
\sm

- We further define the {\em subgradient norm} of $f$ at $x$ by
 \begin{equation} \label{nablamoins}
|\nabla^- f|(x) \: = \: \limsup_{y\to x} \frac{[f(y)-f(x)]_-}{d(x,y)} \: = \:
\limsup_{y\to x} \frac{[f(x)-f(y)]_+}{d(x,y)}.
\end{equation} 
Here $a_+ \: = \: \max(a, 0)$ and $a_- \: = \: \max(-a, 0)$.
Clearly $|\nabla^- f|(x) \: \le \: |\nabla f|(x)$, so the
subgradient norm is a (slightly) finer notion than the gradient norm.
Note that $|\nabla^- f|(x)$ is automatically zero if $f$ has a 
local minimum at $x$. In a sense, $|\nabla^- f|(x)$ measures the 
downward pointing component of $f$ near $x$.
\sm

- Given two probability measures $\mu_0$ and $\mu_1$ on $X$,
the Wasserstein distance (of order~2) $W_2(\mu_0,\mu_1)$ between 
$\mu_0$ and $\mu_1$ is the square root
of the optimal transport cost between $\mu_0$ and $\mu_1$,
when the infinitesimal cost is the square of the distance;
see for instance~\cite[Theorem~7.3]{Villani}.
\sm

- The metric-measure space $(X,d,\nu)$ satisfies a
{\em doubling condition} if the measure $\nu$ is doubling
in the sense of~\cite[eq. (0.1)]{Cheeger (1999)}.
\sm

- The metric-measure space $(X,d,\nu)$ satisfies a
{\em local Poincar\'e inequality} if the measure $\nu$ satisfies
the weak Poincar\'e inequality of type (1,1)
as in~\cite[eq. (4.3)]{Cheeger (1999)}.
\sm

- The metric-measure space $(X,d,\nu)$ is {\em nonbranching}
if any two constant-speed geodesics $[0,1]\to X$ that coincide on an
interval $(t_0,t_1)\subset [0,1]$ are equal.
\med

We will focus on the following three functional inequalities:
\sm

\bul If $K>0$, we say that $(X, d, \nu)$ satisfies a
{\em log Sobolev inequality} with constant $K$, $\LSI(K)$,
if for any $f \in \Lip(X)$ with $\int_X f^2 \: d\nu \: = \: 1$ we have
\begin{equation} \label{ls}
\int_X f^2 \: \log(f^2) \: d\nu \: \le \: \frac{2}{K} \:
\int_X |\nabla^- f|^2 \: d\nu.
\end{equation}
\sm

\bul We say that $(X, d, \nu)$ satisfies a {\em Talagrand inequality} with
constant $K$, $T(K)$, if for any $F \in L^2(X, \nu)$ with 
$\int_X F^2 \: d\nu \: = \: 1$,
we have 
\begin{equation} \label{TI}
W_2(F^2 \nu,\nu) \leq \sqrt{\dps\frac{2 \int_X F^2 \log(F^2) \: d\nu}{K}}.
\end{equation}
\sm

\bul We say that $(X, d, \nu)$ satisfies a (global)
{\em Poincar\'e inequality} with constant $K$, $P(K)$, if for 
any $h \in \Lip(X)$ with $\int_X h\,d\nu = 0$, we have
\begin{equation} \label{GPI}
\int_X h^2\,d\nu \leq 
\frac1{K} \int_X |\nabla^- h|^2\,d\nu. 
\end{equation}

\begin{remark} In~\eqref{ls} and~\eqref{GPI}
we use the subgradient norm defined in~\eqref{nablamoins},
instead of the gradient norm defined in~\eqref{nablaf}. Accordingly, our
log Sobolev and Poincar\'e inequalities are slightly 
stronger statements than those discussed by many other authors.
\end{remark}

Inequalities~\eqref{ls}, \eqref{TI} and~\eqref{GPI} 
are associated with concentration of 
measure~\cite{Toulouse (2000),Bobkov-Gentil-Ledoux 
(2001),Bobkov-Gotze (1999),Ledoux (1999), Ledoux (2001), Villani (StF)}. 
For example, $T(K)$ implies a Gaussian-type concentration of measure.

The following chain of implications, none of which is an equivalence,
is well-known in the context of smooth Riemannian manifolds :
\begin{equation} \label{chain}
[\Ric \geq K ] \Longrightarrow \LSI(K)
\Longrightarrow T(K) \Longrightarrow P(K). 
\end{equation}
A complete proof of~\eqref{chain} is available for instance
in~\cite[Chapters~21 and~22]{Villani (StF)}.
\sm

Our main result is as follows:

\begin{theorem} \label{LSIT}
Let $(X,d,\nu)$ be a compact measured length space.
\sm

(i) If $(X,d,\nu)$ satisfies $T(K)$, for some $K>0$, then it also satisfies $P(K)$.
\sm

(ii) Suppose that $(X,d,\nu)$ satisfies a doubling
condition on the measure and a local Poincar\'e inequality.
If $(X,d,\nu)$ satisfies $\LSI(K)$ for some $K > 0$,
then it also satisfies $T(K)$.
\end{theorem}

It is standard that $\LSI(K)$ implies $P(K)$; see, for example,
\cite[Theorem 6.18]{Lott-Villani}.

The assumptions of Theorem \ref{LSIT}(ii) are satisfied if $(X, d, \nu)$ is
nonbranching and has
Ricci curvature bounded below in the sense of Lott--Villani and Sturm
\cite{Lott-Villani 2,von Renesse,Sturm2}. They are also satisfied if
$(X,d)$ is a length space with Alexandrov curvature bounded
below and Hausdorff dimension $n < \infty$, and $\nu$ is the 
$n$-dimensional Hausdorff measure on $X$; the doubling property follows 
from the Bishop--Gromov 
inequality~\cite[Theorem~10.6.6]{Burago-Burago-Ivanov (2001)} and
the local Poincar\'e inequality was proven in
\cite[Theorem~7.2]{Kuwai-Machigashira-Shioya}.

Theorem~\ref{LSIT} will be proven in Section~\ref{secLS}.
An important technical tool in the proof is the quadratic Hamilton--Jacobi
semigroup, which will be introduced and studied in Section~\ref{secHJ}.
We thank Juha Heinonen for some helpful comments.

\section{Hamilton--Jacobi semigroup} \label{secHJ}

First, we recall the Hamilton-Jacobi semigroup in the case of Riemannian manifolds.
If $M$ is a compact Riemannian manifold, then the quadratic Hamilton--Jacobi 
equation on $M$ is
\begin{equation}
\frac{\partial F}{\partial t} \: + \: \frac{|\nabla F|^2}{2} \: = \: 0. 
\end{equation}
Given an initial condition $f \in C(M)$, the viscosity
solution to the Hamilton--Jacobi equation is given by the Hopf--Lax formula
\begin{equation}\label{HL}
F(t,x) \: = \: \inf_{y\in X} \left [ f(y) + \frac{d(x,y)^2}{2t} \right ],
\end{equation}
where $d$ is the geodesic distance on $M$.
The map that sends $f$ to $F(t, \cdot)$ defines a semigroup
action of $\R_+$ on $C(M)$, called the Hamilton--Jacobi semigroup.

Equation \eqref{HL} does not require any smoothness assumption,
so the following definition makes sense.

\begin{definition}
Let $(X,d)$ be a compact metric space. Given $f \in C(X)$ and 
$t\geq 0$, we define a map $Q_t \: : \: X \rightarrow \R$ by
\begin{equation} \label{Qt}
(Q_tf)(x) = \inf_{y\in X} \left [ f(y) + \frac{d(x,y)^2}{2t} \right ],
\end{equation}
with the convention that $Q_0 f = f$.
\end{definition}

If $X$ is a length space, then the map $Q_t$ defines a semigroup
action of $\R^+$ on $C(X)$; see part (i) of Theorem \ref{propHJ} below.
We may then speak of the ``Hamilton--Jacobi semigroup''.
The next theorem establishes some of its basic properties.

\begin{theorem} \label{propHJ} 

(i) For any $s, t \ge 0$,  $Q_t Q_s f \: = \: Q_{t+s}f$.

(ii) For any $x\in X$, $\inf f \leq (Q_tf)(x) \leq f(x)$.

(iii) For any $t>0$, $Q_tf \in \Lip(X)$.

(iv) For any $x\in X$, $(Q_tf)(x)$ is a nonincreasing function of $t$,
that converges monotonically to $f(x)$ as $t\to 0$. In particular, 
$\lim_{t \rightarrow 0} Q_tf \: = \: f$ in $C(X)$.

(v) For any $t\geq 0$, $s>0$ and $x\in X$,
\begin{equation}
\frac{|Q_{t+s}f(x) - Q_tf(x)|}{s} \leq \frac{\|Q_tf\|_{\Lip}^2}{2}.
\end{equation}

(vi) For any $x\in X$ and $t \geq 0$,
\begin{equation} \label{liminfeqn}
\liminf_{s \rightarrow 0^+} \frac{(Q_{t+s}f)(x)-(Q_tf)(x)}s \ge
- \frac{|\nabla^- Q_tf|(x)^2}2.
\end{equation}

(vii) If $(X,d,\nu)$ satisfies a doubling
condition on the measure and a local Poincar\'e inequality, then
for $t > 0$ and $\nu$-almost any $x\in X$,
\begin{equation} 
\lim_{s \rightarrow 0^+} \frac{(Q_{t+s}f)(x)-(Q_tf)(x)}s =
- \frac{|\nabla^- Q_tf|^2(x)}2.
\end{equation}

(viii) If $(X,d)$ is a finite-dimensional space with 
Alexandrov curvature bounded below then
for any $t > 0$ and any $x\in X$,
\begin{equation} 
\lim_{s \rightarrow 0^+} \frac{(Q_{t+s}f)(x)-(Q_tf)(x)}s =
- \frac{|\nabla^- Q_tf|^2(x)}2.
\end{equation}
\end{theorem}

\begin{remark} 
Part (vii) of
Theorem \ref{propHJ} will be used in the proof of Theorem~\ref{LSIT}.
Part (viii) is not needed for the proof of Theorem \ref{LSIT},
but may be of independent interest.  Parts (vii) and (viii)
show that for $t > 0$, the function 
$F(t,x)=Q_t f(x)$
 satisfies the Hamilton-Jacobi equation
\begin{equation}
\frac{\partial F}{\partial t} \: + \: \frac{|\nabla^- F|^2}{2} \: = \: 0,
\end{equation}
almost everywhere in the case of (vii) and everywhere in the case of (viii).
\end{remark}

\begin{remark} Theorem~\ref{propHJ} is reminiscent of known properties
of Hamilton--Jacobi equations in a smooth setting; see
e.g.~\cite{Cannarsa-Sinestrari (2004)}.
However, even in the context of Riemannian manifolds,
we have been unable to find exactly this statement in the literature.
On the one hand, the vast majority of works are only concerned with Euclidean 
or Hilbert spaces. On the other hand, the use of the subgradient norm 
is a bit nonstandard.
\end{remark}

\begin{remark} More general Hamilton--Jacobi semigroups will be considered
in~\cite[Appendix of Chapter~22]{Villani (StF)}, of the form
\[ Q_tf(x) = \inf_{y\in X}\: \left[ f(y) + 
t\, L \left(\frac{d(x,y)}{t}\right)\right],\]
where $L:\R_+\to\R_+$ is increasing, convex and locally semiconcave, with $L(0)=0$.
Theorem~\ref{propHJ} can be extended {\em mutatis mutandis}
to this more general situation (apart maybe from (viii)).
In this generalization, a few minor complications arise 
if $L'(\infty)<+\infty$. (For simplicity,
only Riemannian manifolds are considered in~\cite[Chapter~22]{Villani (StF)},
but nonsmooth spaces can be treated as in the present paper.)
At the level of geometric applications, an interesting case
in relation to Poincar\'e inequalities (as opposed to
log Sobolev inequalities) is when $L(s)$ is asymptotic to $s^2$ for
small $s$, and to $s$ for large $s$.
\end{remark}

\begin{proof}[Proof of Theorem \ref{propHJ}]
To prove (i), we first claim that for all $x, y \in X$ and $s,t>0$,
\begin{equation} 
\label{fmlinfinf} 
\frac{d(x,y)^2}{t+s} \: = \:
\inf_{z \in X} \left [ \frac{d(x,z)^2}t + \frac{d(z,y)^2}s \right ].
\end{equation}
The triangle inequality implies that the left-hand side 
of~\eqref{fmlinfinf} is less than or equal to the right-hand side.
The equality in (\ref{fmlinfinf}) comes from choosing a minimal
geodesic between $x$ and $y$, and a point $z$ on this
geodesic with $d(x,z) \: = \: \frac{t}{s+t} \: d(x,y)$.

{From} (\ref{fmlinfinf}), we obtain
\begin{equation}
(Q_{t+s} f)(x) =  
\inf_{y\in X} \left [ f(y) + \frac{d(x,y)^2}{2(t+s)} \right ]
= \inf_{y\in X} \inf_{z\in X} \left [ f(y) + \frac{d(x,z)^2}{2t}
+ \frac{d(z,y)^2}{2s} \right ] = (Q_t Q_s f)(x), 
\end{equation}
which proves (i).
\sm

For part (ii), the inequality on the left is obvious, while the
inequality on the right follows from the choice $y=x$ in the definition
of $(Q_tf)(x)$.
\sm

Part (iii) follows from
\begin{align}
(Q_tf)(x) - (Q_tf)(x') & \leq \frac1{2t} \sup_{y\in X}
[ d(x,y)^2-d(x',y)^2] \\
& \leq \left ( \frac1{2t} \sup_{y\in X} 
[d(x,y)+d(x',y)] \right ) d(x,x') \notag \\
& \leq \frac{\diam(X)}{t} \: d(x,x').\notag
\end{align}

In view of (i) and (ii), for any $s,t>0$ and $x \in X$,
\begin{equation} 
(Q_{t+s} f)(x) \leq (Q_tf)(x),
\end{equation}
so $(Q_tf)(x)$ is indeed a nonincreasing function of $t$.
Given $f \in C(X)$, put $C=C(f)=2(\sup f-\inf f)$.
If $y$ is such that $d(x,y) \ge \sqrt{Ct}$ then
\begin{equation} 
f(y) + \frac{d(x,y)^2}{2t} \ge (\inf f) + \frac{C}{2}
= \sup f \geq f(x).
\end{equation}
We conclude that
\begin{equation} 
(Q_tf)(x) = \inf_{y \in B_{\sqrt{Ct}}(x)} 
\left [ f(y) + \frac{d(x,y)^2}{2t}
\right ].
\end{equation}

Given $x \in X$ and $\var>0$, choose $\delta > 0$ so that
\begin{equation} 
d(x,y) <\delta \Longrightarrow |f(x)-f(y)|<\var.
\end{equation}
If $t \: \leq \: \frac{\delta^2}{C}$ then $\sqrt{Ct} \: \le \: \delta$, so
\begin{equation} 
(Q_tf)(x) \geq \inf_{y \in B_{\delta}(x)} 
\left [ f(y) + \frac{d(x,y)^2}{2t}
\right ] \geq f(x)-\var. 
\end{equation}
This shows that $\lim_{t \rightarrow 0} (Q_tf)(x) \: = \: f(x)$.
Since the convergence is monotone and $X$ is compact,
the convergence is uniform. This proves part (iv)
of the theorem.
\sm

Next, for $g \in C(X)$, 
we write (with $C = C(g)$ and the convention that $0\cdot\infty=0$)
\begin{align*} \frac{g(x) - (Q_sg)(x)}s & = 
\frac1{s} \sup_{y \in B_{\sqrt{Cs}}(x)} \left [ g(x) - g(y) - 
\frac{d(x,y)^2}{2s} \right ] \\
& \leq \sup_{y \in B_{\sqrt{Cs}}(x)}
\left ( \frac{[g(x)-g(y)]_+}{d(x,y)}\,\frac{d(x,y)}s
- \frac{d(x,y)^2}{2s^2} \right ) \\
& \leq \sup_{y \in B_{\sqrt{Cs}}(x)} 
\frac12 \left ( \frac{[g(x)-g(y)]_+}{d(x,y)} \right )^2.
\end{align*}
If $g=Q_tf$, in view of (i) and (ii) this becomes
\begeq\label{ineqQtf}
0 \leq \frac{Q_tf(x) - (Q_{t+s}f)(x)}s \leq \sup_{y \in B_{\sqrt{Cs}}(x)} 
\frac12 \left ( \frac{[Q_tf(x)-Q_tf(y)]_+}{d(x,y)} \right )^2.
\endeq
Then statement (v) follows immediately. If now we let $s\to 0^+$ then
the definition of $|\nabla^-Q_tf|$ implies that
\begeq \label{limsupQt} 
\limsup_{s \rightarrow 0^+} \sup_{y \in B_{\sqrt{Cs}}(x)} 
\frac12 \left ( \frac{[Q_tf(x)-Q_tf(y)]_+}{d(x,y)} \right )^2
\leq \frac{|\nabla^- Q_tf|^2(x)}2, 
\endeq
and (vi) is also true.
\sm

We now turn to (vii) and (viii), which are the most delicate parts of
the theorem. Again with $g=Q_tf$, we want to prove that
\begin{equation} \label{needed}
\liminf_{s \rightarrow 0^+} \left [ \frac{g(x) - (Q_{s}g)(x)}s \right ]
\geq \frac{|\nabla^-g|^2(x)}{2}. 
\end{equation}
The case when $|\nabla^-g|(x)=0$ is obvious, since $(Q_{t}g)(x)$ is
a nonincreasing function of $t$. So in what follows we assume that
$|\nabla^-g|(x)>0$. 

We write
\begin{align} \label{mess}
\frac{g(x) - (Q_{s}g)(x)}s & = 
\frac1{s} \sup_{y \in X} \left [ g(x) - g(y) - \frac{d(x,y)^2}{2s} \right ] \\
& \geq \sup_{y \in S_{s|\nabla^- g|(x)}(x)}
\left ( \left [ \frac{g(x)-g(y)}{d(x,y)} \right ]
|\nabla^-g|(x) - \frac{|\nabla^-g|^2(x)}2 \right ). \notag
\end{align}
Put
\begin{equation}
\psi(r) = \sup_{y \: \in S_r(x)}
\frac{g(x)-g(y)}{d(x,y)}.
\end{equation}
As $\limsup_{r \rightarrow 0^+} \psi(r) \: = \: |\nabla^- g|(x) \: > \: 0$,
if we can show that 
$\liminf_{r \rightarrow 0^+} \psi(r) \: = \: |\nabla^- g|(x)$ then
equation (\ref{mess}) will imply (\ref{needed}).

For (vii), we use results from \cite{Cheeger (1999)}. 
By~\cite[Theorem 10.2]{Cheeger (1999)}, 
the Lipschitz function $g$ admits generalized linear derivatives
$g_{0,x}$ at $x$ for 
$\nu$-almost
all $x \in X$. For such an $x$, suppose that there is a 
sequence $r_i \rightarrow 0$ such that $\lim_{i \rightarrow 0}
\psi(r_i) \: = \: |\nabla^- g|(x) - \epsilon$ for some $\epsilon > 0$.
After passing to a further subsequence, we can assume that
the rescaled measured length spaces $(X, x, r_i^{-1/2} d, \nu)$
converge to a measured tangent cone $(X_x, x_\infty, d_\infty, \nu_\infty)$
and the rescaled functions 
$g_{r_i^{1/2}, x} = (g-g(x))/r_i^{1/2}$
converge to a generalized linear function $g_{0,x}$ on $X_\infty$; 
see \cite[Theorem 10.2]{Cheeger (1999)}.
(Note that we rescale by $r_i^{1/2}$ and not $r_i$; any
rate $s_i\to 0$ such that $r_i=o(s_i)$ would do.)
Then $|\nabla^- g_{0,x}|(x_\infty) \: = \: |\nabla^- g|(x) - \epsilon$. 
From \cite[Theorem 8.10]{Cheeger (1999)}, there is a unit-speed line
$\gamma$ in $X_\infty$ through $x_\infty$ which is an integral curve for
$g_{0,x}$. That is, $\gamma(0) \: = \: x_\infty$ and
$\frac{d}{dt} g_{0,x}(\gamma(t)) \: = \: |\nabla g_{0,x}|$. It follows that
$|\nabla^- g_{0,x}|(x_\infty) \: \geq \: |\nabla g_{0,x}|(x_\infty)$. 
However,
from \cite[Theorem 10.2]{Cheeger (1999)}, one has
$|\nabla g_{0,x}|(x_\infty) \: = \: |\nabla g|(x)$. Thus 
$|\nabla g|(x) \: \leq \: |\nabla^- g|(x) - \epsilon$, 
which is a contradiction. This proves statement (vii).

\begin{remark} This reasoning shows actually shows that $|\nabla^-g|(x)\geq
|\nabla g|(x)$, so $|\nabla^-g|(x)=|\nabla g|(x)$ (for 
$\nu$-almost all $x$).
\end{remark}

Statement (viii), with convergence for all $x \in X$, requires additional
regularity for $X$. We will use the notion of {\em quasigeodesics} in 
Alexandrov spaces, as studied in \cite{Perelman-Petrunin}.
The following properties will be useful: 
(a) squared distance functions, when restricted to quasigeodesics,
satisfy the same curvature-dependent differential inequalities as
when restricted to geodesics (inequality~\eqref{ineqgeod} below);
(b) quasigeodesics can be extended to all positive times;
(c) uniform limits of quasigeodesics are quasigeodesics, and this
statement goes through for quasigeodesics defined on a 
Gromov--Hausdorff converging sequence of Alexandrov spaces.
We recall that by definition, nontrivial quasi-geodesics are
parametrized by arc-length. 

\begin{lemma} \label{quasi}
Let $X$ be a finite-dimensional compact length space with Alexandrov 
curvature bounded below. Fix $x \in X$. Then
\sm

(i) There is some $\delta > 0$ so that
each complete quasigeodesic $\gamma \: : \: [0, \infty) \rightarrow X$
starting from $x$ intersects $S_\delta(x)$. 
\sm

(ii) There is a function $\sigma \: : \: (0, \delta) \rightarrow \R_+$ with
$\lim_{r \rightarrow 0^+} \sigma(r) \: = \: 0$ so that if 
$\gamma \: : \: [0,L] \rightarrow X$ is a quasigeodesic segment
starting from $x$ with $\gamma(L) \in S_r(x)$ and
$\gamma([0,L]) \subset \overline{B_r(x)}$ then
$\left| \frac{L}{r} - 1 \right|
\: \le \: \sigma(r)$.
\end{lemma}

\begin{remark} Of course, if $X$ is a Riemannian manifold then
this lemma holds true for geodesics. (Take $\delta$ to be the
injectivity radius at $x$ and take $\sigma=0$.) 
\end{remark}

\begin{proof}[Proof of Lemma~\ref{quasi}]
Suppose that (i) is not true.  Then for each $i \in \Z^+$, there is a 
quasigeodesic $\gamma_i \: : \: [0, \infty) \rightarrow X$ starting
from $x$ that remains in $B_{1/i}(x)$. Taking a convergent
subsequence of $\{\gamma_i\}_{i=1}^\infty$
\cite[\S 2]{Perelman-Petrunin} gives a quasigeodesic 
$\gamma_\infty \: : \: [0, \infty) \rightarrow X$ whose image is $\{x\}$.
This contradicts the fact that a quasigeodesic has unit speed.

Suppose that (ii) is not true.  Then there is an $\epsilon > 0$ along with 
a sequence $\{r_i\}_{i=1}^\infty$ converging to zero and a sequence
of quasigeodesic segments $\gamma_i \: : \: [0, L_i] \rightarrow X$
starting from $x$ so that for all $i \in \Z^+$, we have that
$\gamma_i(L_i) \in
S_{r_i}(x)$, $\gamma_i([0, L_i]) \subset \overline{B_{r_i}(x)}$
and $\frac{L_i}{r_i} \: \ge \: 1+ \epsilon$.
Rescaling the pointed Alexandrov space $(X, x)$ by $\frac{1}{r_i}$ and
taking a convergent
subsequence of $\{\gamma_i\}_{i=1}^\infty$
\cite[Theorem 2.2]{Perelman-Petrunin}, we obtain a
quasigeodesic segment $\gamma_\infty \: : \: [0, L_\infty] \rightarrow
C_xX$ starting at the vertex $o$ of the tangent cone $C_x X$ so that
$\gamma_\infty(L_\infty) \in S_1(o)$ (if $L_\infty < \infty)$,
$\gamma_\infty([0,L_\infty]) \subset \overline{B_1(o)}$ and
$L_\infty \ge 1+\epsilon$. However, one can check that a
quasigeodesic in $C_xX$ starting at $o$ must be a radial geodesic,
which is a contradiction.
\end{proof}

Let us go back to the proof of Theorem~\ref{propHJ}, part (viii).
As $X$ is compact with Alexandrov curvature bounded below,
there is a $K \ge 0$ so that
for all quasigeodesic segments $\gamma \: : \: [0,u] \rightarrow X$
starting from $x$, all $u^\prime \in [0, u]$ 
and all $z \in X$,
\begin{equation} \label{ineqgeod}
d(\gamma(u^\prime), z)^2 \: - \: \frac{u^\prime}{u} \: 
d(\gamma(u), z)^2 \: - \:
\frac{u-u^\prime}{u} \: d(x, z)^2 \: \ge \: - \: K \: u^\prime (u-u^\prime).
\end{equation}
Then
\begin{equation}
\frac{1}{2t} \: \inf_{z \in X} \left( 
d(\gamma(u^\prime), z)^2 \: - \: \frac{u^\prime}{u} \: 
d(\gamma(u), z)^2 \: - \:
\frac{u-u^\prime}{u} \: d(x, z)^2 \right) \: \ge \: - \:
\frac{1}{2t} \:  K \: u^\prime (u-u^\prime).
\end{equation}
As
\begin{align}
g(\gamma(u^\prime)) \: - \: \frac{u^\prime}{u} \: g(\gamma(u)) \: - \: 
\frac{u-u^\prime}{u} \: g(x) \: = \: &
\inf_{z^\prime \in X} \sup_{z,w \in X}
\left[ f(z^\prime) \: - \: \frac{u^\prime}{u} \: f(z) \: - \: 
\frac{u-u^\prime}{u} \: f(w) \: + \right. \\
& \left. 
\frac{d(\gamma(u^\prime), z^\prime)^2}{2t} \: - \:
\frac{u^\prime}{u} \: \frac{d(\gamma(u), z)^2}{2t} \: - \:
\frac{u-u^\prime}{u} \: \frac{d(x, w)^2}{2t} \right], \notag
\end{align}
by considering the case when $z^\prime \: = \: z \: = \: w$, we obtain
\begin{equation}
g(\gamma(u^\prime)) \: - \: \frac{u^\prime}{u} \: g(\gamma(u)) \: - \: 
\frac{u-u^\prime}{u} \: g(x) \: \ge \: - \:
\frac{1}{2t} \:  K \: u^\prime(u-u^\prime).
\end{equation}
Equivalently,
\begin{equation} \label{equivalently}
\frac{g(x)-g(\gamma(u))}{u} \: + \: \frac{1}{2t} \: K \: u \: \ge \:
\frac{g(x)-g(\gamma(u^\prime))}{u^\prime} \: + \: \frac{1}{2t} \: K \:
u^\prime.
\end{equation}

In  order to prove that $\liminf_{r \rightarrow 0^+} \psi(r) \: = \: |\nabla^- g|(x)$,
suppose that $\liminf_{r \rightarrow 0^+} \psi(r) \: = \: |\nabla^- g|(x)
\: - \: \epsilon$
for some $\epsilon > 0$. Then there are sequences 
$\{u_i^\prime\}_{i=1}^\infty$ and $\{v_j\}_{j=1}^\infty$
converging to zero with
\begin{equation} \label{firstlimit}
\lim_{i \rightarrow \infty} \psi(u_i^\prime) \: = \:  |\nabla^- g|(x)
\end{equation}
and
\begin{equation} \label{secondlimit}
\lim_{i \rightarrow \infty} \psi(v_i) \: = \:  |\nabla^- g|(x) -
\epsilon.
\end{equation}
We may assume that $u_i^\prime < v_i < \delta$, where
$\delta$ is from Lemma \ref{quasi}(i).

In particular, there are points
$y_i^\prime \in S_{u_i^\prime}(x)$ so that
\begin{equation} \label{thirdlimit}
\lim_{i \rightarrow \infty} \frac{g(x)-g(y_i^\prime)}{u_i} \: = \: 
|\nabla^- g|(x).
\end{equation}
Choose a minimizing geodesic $\gamma_i$ from
$x$ to $y_i^\prime$. Extend it to a complete quasigeodesic
$\gamma_i \: : \: [0, \infty) \rightarrow X$. Put
\begin{equation}
u_i \: = \: \inf \{ w_i \: : \: \gamma_i(w_i) \in S_{v_i}(x) \}.
\end{equation}
{From} Lemma \ref{quasi}(i), $u_i$ exists.
As $\gamma_i$ is parametrized by arclength,
$u_i \ge v_i > u_i^\prime$.
{From} (\ref{equivalently}),
\begin{equation}
\frac{v_i}{u_i} \:
\frac{g(x)-g(\gamma(u_i))}{v_i} \: + \: \frac{1}{2t} \: K \: u_i \: \ge \:
\frac{g(x)-g(y_i^\prime)}{u_i^\prime} \: + \: \frac{1}{2t} \: K \:
u_i^\prime.
\end{equation}
In particular,
\begin{equation} \label{particular}
\frac{v_i}{u_i} \:
\psi(v_i) \: + \: \frac{1}{2t} \: K \: u_i \: \ge \:
\frac{g(x)-g(y_i^\prime)}{u_i^\prime} \: + \: \frac{1}{2t} \: K \:
u_i^\prime.
\end{equation}
{From} Lemma \ref{quasi}(ii),
\begin{equation}
\lim_{i \rightarrow \infty} \frac{v_i}{u_i} \: = \: 1.
\end{equation}
Taking $i \rightarrow \infty$ in (\ref{particular}), 
we get a contradiction to (\ref{secondlimit}) and (\ref{thirdlimit}).
\end{proof}

\section{Proof of Theorem \ref{LSIT}} \label{secLS}

Armed with Theorem~\ref{propHJ}, we can now use the strategy
of~\cite{Bobkov-Gentil-Ledoux (2001)} to prove Theorem~\ref{LSIT}.

\begin{proof}[Proof of Theorem~\ref{LSIT}, part (i)]
Let $h \in \Lip(X)$ satisfy $\int_X h\,d\nu=0$. Introduce
\begin{equation} 
\psi(t) = \int_X e^{K t Q_t h}\,d\nu.
\end{equation}
{From} Talagrand's inequality in its dual formulation
(see \cite[p. 16]{Bobkov-Gotze (1999)},
\cite[Exercise 9.15]{Villani} or \cite[Chapter~22]{Villani (StF)}), 
we know that $\psi(t) \leq \exp(Kt \int_X h\,d\nu)=1$. Hence $\psi$ has
a maximum at $t=0$. Combining this with $\int h\,d\nu =0$, we find
\begeq\label{limsuppsi}
0 \leq \limsup_{t\to 0^+} \left(\frac{1-\psi(t)}{Kt^2}\right)
= \limsup_{t\to 0^+} \int_X \left( \frac{1 + Kt\,h - e^{KtQ_th}}{Kt^2}\right)\,d\nu.
\endeq

By the boundedness of $Q_th$ and Theorem~\ref{propHJ}(iv), 
\begin{align}
e^{Kt Q_th} & = 1 + Kt Q_th + \frac{K^2 t^2}2\, (Q_th)^2 + O(t^3) \\
& = 1 + K t Q_th + \frac{K^2 t^2}2\, h^2 + o(t^2). \notag
\end{align}
So the right-hand side of~\eqref{limsuppsi} equals
\begin{equation}
\limsup_{t\to 0^+} \int_X \left( \frac{h-Q_th}{t} \right)\,d\nu
- \frac{K}2 \int_X h^2\,d\nu.
\end{equation}
By Theorem~\ref{propHJ}(v), $(h-Q_th)/t$ is bounded, which allows us to
apply Fatou's lemma in the form
\begin{equation}
\limsup_{t\to 0^+} \int_X \left( \frac{h-Q_th}{t} \right)\,d\nu
\leq \int_X \limsup_{t\to 0^+} \left( \frac{h-Q_th}{t} \right)\,d\nu. 
\end{equation}
Then Theorem~\ref{propHJ}(vi) implies that
\begin{equation}
 \int_X \limsup_{t\to 0^+} \left( \frac{h-Q_th}{t} \right)\,d\nu
\leq \int_X \frac{|\nabla^- h|^2}2\,d\nu.
\end{equation}
All in all, the right-hand side of~\eqref{limsuppsi} can be bounded above by
\begin{equation}
\frac12 \int_X |\nabla^- h|^2\,d\nu - \frac{K}2 \int_X h^2\,d\nu,
\end{equation}
so this expression is nonnegative. This concludes the proof.
\end{proof}

\begin{proof}[Proof of Theorem~\ref{LSIT}, part (ii)]
{From} Talagrand's inequality in its dual formulation,
it is sufficient
to show that for all $g\in C(X)$,
\begin{equation} 
\int_X e^{K \inf_y [ g(y) + \frac{d(x,y)^2}2 ]}\,d\nu(x) \leq
e^{K\int_X g\,d\nu}.
\end{equation}

Put
\begin{equation} 
\phi(t)= \frac{1}{K t} \log \left ( \int_X e^{K t Q_t g}
\,d\nu \right ). 
\end{equation}
Since $g$ is bounded, Theorem~\ref{propHJ}(ii) implies that $Q_tg$ is
bounded, uniformly in $t$. Thus
\begin{equation}
\int_X e^{KtQ_tg} \,d\nu   = 1 + Kt \int_X Q_tg\,d\nu  + O(t^2)
\end{equation}
and
\begin{equation} 
\phi(t) = \int_X Q_tg\,d\nu + O(t).
\end{equation}
By Theorem~\ref{propHJ}(iv), $Q_tg$ converges uniformly to $g$ as 
$t\to 0^+$, and so
\begin{equation} 
\lim_{t \rightarrow 0^+} \phi(t) \: = \: \int_X g\,d\nu. 
\end{equation}
Therefore, our goal will be achieved if we can show that
$\phi(1) \leq \lim_{t\to 0^+} \phi(t)$.
For this, it suffices to show that $\phi(t)$ is nonincreasing in $t$.

Let $t\in (0,1]$ be given. For $s>0$, we have 
\begin{align} \label{2terms}
\frac{\phi(t+s)-\phi(t)}s \: = \: &
\frac1{s} \left( \frac1{K(t+s)} - \frac1{Kt} \right)
\log \int_X e^{K(t+s)Q_{t+s}g}\,d\nu \: + \\
& \frac1{Kts} \left( \log \int_X e^{K(t+s)Q_{t+s}g}\,d\nu
- \log\int_X e^{KtQ_{t}g}\,d\nu \right). \notag
\end{align}
As $s\to 0^+$, $e^{K(t+s) Q_{t+s}g}$ converges uniformly to 
$e^{KtQ_tg}$. Thus the limit of the first term in the right-hand side above,
as $s\to 0^+$, is
\begin{equation}
- \: \frac1{K t^2} \: \log
\left ( \int_X e^{K t Q_tg}\,d\nu \right ),
\end{equation}
while the limit of the second term is
\begeq\label{secondterm} 
 \frac1{Kt\int e^{KtQ_{t}g}\,d\nu} \;
\lim_{s\to 0^+}\: \left[\frac1{s} \left( \int_X
e^{K(t+s)Q_{t+s}g}\,d\nu
- \int_X e^{Kt Q_{t}g}\,d\nu\right) \right],
\endeq
provided that the latter limit exists. We rewrite the expression
inside the square brackets as
\begeq\label{2terms2}
\int_X \left(\frac{e^{K(t+s)Q_{t+s}g} - e^{Kt Q_{t+s}g}}s\right)
\,d\nu
+ \int_X  \left( \frac{e^{KtQ_{t+s}g} - e^{Kt Q_{t}g}}s\right)\,d\nu.
\endeq

The integrand of the first term in~\eqref{2terms2} can be rewritten as
$(e^{Kt Q_{t+s}g}) ( e^{KsQ_{t+s}g} -1)/s$,
which converges uniformly to $(e^{Kt Q_tg}) K Q_tg$ as $s\to 0^+$.
So the first integral in~\eqref{2terms2} converges to
$\int_X (K Q_tg) e^{K tQ_tg}\,d\nu$.

We now turn to the second term of~\eqref{2terms2}.
By Theorem~\ref{propHJ}(vii), for 
$\nu$-almost all $x\in X$ we have
\begin{equation}
Q_{t+s}g(x) = Q_tg(x) - s \left(\frac{|\nabla^-Q_tg(x)|^2}2 + o(1)\right),
\end{equation}
and therefore
\begeq\label{tcd1} 
\lim_{s \rightarrow 0^+}
\frac{e^{KtQ_{t+s}g(x)} - e^{KtQ_tg(x)}}{s} \: = \:
- \: Kte^{KtQ_tg} \frac{|\nabla^- Q_tg(x)|^2}2. 
\endeq
On the other hand, parts (iv) and (v) of
Theorem~\ref{propHJ} imply that
\begin{equation} Q_{t+s} g = Q_tg + O(s).
\end{equation}
Since $Q_tg(x)$ is uniformly bounded in $t$ and $x$, we deduce that
\begeq\label{tcd2} \frac{e^{Kt Q_{t+s}g} - e^{KtQ_tg}}{s} = O(1) 
\endeq
as $s\to 0^+$. The combination of~\eqref{tcd1} and \eqref{tcd2} makes
it possible to pass to the limit by dominated convergence, to obtain
\begin{equation} 
\lim_{s \rightarrow 0^+}
\int_X  \left( \frac{e^{KtQ_{t+s}g} - e^{Kt Q_{t}g}}s\right)
\,d\nu  \: = \: - \: 
Kt\int_X \frac{|\nabla^- Q_tg|^2}2\,e^{K t Q_tg}\,d\nu.
\end{equation}

In summary,
\begin{align}
\lim_{s\to 0^+} \left[\frac{\phi(t+s)-\phi(t)}s\right]
\: = \: & \frac1{K t^2 \int_X e^{K t Q_tg}\,d\nu}\left [ - \: 
\left ( \int_X e^{K t Q_tg}\,d\nu \right )\log
\left ( \int_X e^{K t Q_tg}\,d\nu \right ) \: + \: \right. \\
& \left. \int_X (K t Q_tg) e^{K tQ_tg}\,d\nu \: - \:
\frac{1}{2K}\int_X \bigl(Kt|\nabla^- Q_tg|\bigr)^2\,e^{K t Q_tg}\,d\nu
\right ]. \notag
\end{align}
Inequality $\LSI(K)$ implies that this quantity is nonpositive, 
which concludes the proof.
\end{proof}


\begin{thebibliography}{10}

\bibitem{Toulouse (2000)}
C. An\'e, S. Blach\`ere, D.~Chafa\"{\i}, P.~Foug\`eres, I.~Gentil,
F.~Malrieu, C.~Roberto and G.~Scheffer,
\underline{Sur les in\'egalit\'es de Sobolev logarithmiques},
Panoramas et Synth\`eses 10, Soci\'et\'e Math\'ematique de France
(2000)

\bibitem{Bobkov-Gentil-Ledoux (2001)}
S. Bobkov, I. Gentil and M. Ledoux,
``Hypercontractivity of Hamilton-Jacobi equations'',
J. Math. Pures Appl. 80, p. 669-696 (2001)

\bibitem{Bobkov-Gotze (1999)}
S.G. Bobkov and F. G\"otze,
``Exponential integrability and transportation cost related to
logarithmic Sobolev inequalities'',
J. Funct. Anal. 163, p. 1-28 (1999)

\bibitem{Burago-Burago-Ivanov (2001)}
D. Burago, Y. Burago and S. Ivanov,
\underline{A course in metric geometry},
Graduate Studies in Mathematics 33, American Mathematical Society,
Providence (2001)

\bibitem{Cannarsa-Sinestrari (2004)}
P. Cannarsa and C. Sinestrari,
\underline{Semiconcave functions, Hamilton--Jacobi equations,
and optimal control},
Progress in Nonlinear Differential Equations and Applications 58,
Birkh\"auser, Boston (2004)

\bibitem{Cheeger (1999)} J. Cheeger,
``Differentiability of Lipschitz functions on metric measure spaces'',
Geom. and Funct. Anal. 9, p. 428-517 (1999)

\bibitem{Kuwai-Machigashira-Shioya} K. Kuwai, Y. Machigashira and T. Shioya,
``Sobolev spaces, Laplacian, and heat kernel on Alexandrov spaces'',
Math. Z. 238, p. 269-316 (2001)

\bibitem{Ledoux (1999)}
M. Ledoux, ``Concentration of measure and logarithmic Sobolev inequalities'',
S\'eminaire de Probabilit\'es XXXIII, 
Lecture Notes in Mathematics 1709, p. 120-216
Springer-Verlag, Berlin (1999)

\bibitem{Ledoux (2001)}
M. Ledoux, \underline{The concentration of measure phenomenon},
American Mathematical Society, Providence (2001)

\bibitem{Lott-Villani} J. Lott and C. Villani, 
``Ricci curvature for metric-measure spaces via optimal transport'', 
to appear, Ann. of Math., 
http://www.arxiv.org/abs/math.DG/0412127

\bibitem{Lott-Villani 2} J. Lott and C. Villani, ``Weak curvature conditions and
functional inequalities'', to appear, J. of Funct. Anal., http://www.arxiv.org/abs/math.DG/0506481

\bibitem{Otto-Villani (2000)} F. Otto and C. Villani,
``Generalization of an inequality by Talagrand, and links with the
logarithmic Sobolev inequality'',
J. Funct. Anal. 173, p. 361-400 (2000)

\bibitem{Perelman-Petrunin} G. Perelman and A. Petrunin,
``Quasigeodesics and gradient curves in Alexandrov spaces'',
unpublished preprint

\bibitem{von Renesse} M. von Renesse, ``On local Poincar\'e via transportation'',
to appear, Math. Zeitschrift, http://www.arxiv.org/abs/math.MG/0505588

\bibitem{Sturm1} K.-T Sturm, ``On the geometry of metric measure spaces I '',
Acta Math. 196, p. 65-131 (2006)

\bibitem{Sturm2} K.-T Sturm, ``On the geometry of metric measure spaces II '',
Acta Math. 196, p. 133-177 (2006)

\bibitem{Villani} C. Villani, \underline{Topics in Optimal Transportation},
Graduate Studies in Mathematics 58, American  Mathematical 
Society, Providence (2003)

\bibitem{Villani (StF)} C. Villani, \underline{Optimal Transport,
Old and New}, in preparation

\end{thebibliography}
\end{document}